\documentclass[10pt,amssymb,amscd]{amsart}
\usepackage{color}
\usepackage{amssymb}
\usepackage{amsmath}
\usepackage{amsthm} 
\usepackage{mathrsfs}                                                            
\usepackage[all]{xy}
\usepackage{wrapfig}

\def\2{{1\over 2}}

\newcommand{\rf}[1]{(\ref{#1})}

\newcommand{\p}{\partial}

\def\p{\partial}

\def\b{\bar}

\def\<{\langle}
\def\>{\rangle}
\def\+{\dagger}

\begin{document}
\title{On the continuous series for $\widehat{sl(2,\mathbb{R})}$ }
\author{Igor B. Frenkel}
\address{\newline Department of Mathematics, \newline
Yale University,\newline
10 Hilhouse ave,\newline
New Haven, CT, 06511\newline
igor.frenkel@yale.edu}
\author{Anton M. Zeitlin}
\address{ 
\newline Department of Mathematics,\newline
Columbia University,
\newline 2990 Broadway, New York,\newline
NY 10027, USA.
\newline
zeitlin@math.columbia.edu\newline
http://math.columbia.edu/$\sim$zeitlin \newline
http://www.ipme.ru/zam.html  }

\begin{abstract}
We construct the representations of $\widehat{sl(2,\mathbb{R})}$ 
starting from the unitary representations of the loop $ax+b$-group. Our approach involves a combinatorial 
analysis of the correlation functions of the generators and renormalization of the appearing divergencies. We view our construction as a step towards a realization of the principal series representations of $\widehat{sl(2,\mathbb{R})}$.
\end{abstract}
\maketitle
\section{Introduction}

The principal series representations of split reductive real Lie groups is one of the cornerstones of the classical representation theory. Recently a proper quantum analogue of the principal series in the case of the quantum group $U_q(sl(2,\mathbb{R}))$, or more precisely of its modular double \cite{faddeev},\cite{kls} has been found in \cite{teschner}, \cite{teschner2}, based on the previous work in \cite{schmudgen}. These unitary representations, which do not have a classical limit, but behave similarly to the finite-dimensional representations of the quantum group $U_q(su(2))$, in particular, they form a tensor category. A generalization of these representations to higher rank quantum groups has been also obtained in \cite{gkl},\cite{gerasimov},\cite{fip}, \cite{ivan2}, \cite{ivan3}.  

One of the fundamental results in the theory of quantum groups is the equivalence of braided tensor categories of finite-dimensional representations of quantum groups and of dominant highest weight representations of affine Lie algebras, associated with the same simple Lie algebra \cite{ms}, \cite{KL}.
By Drinfeld-Sokolov reduction, this equivalence can be extended to the so-called W-algebras, which in the simplest 
case is nothing but the Virasoro algebra. It was shown in \cite{teschner} that the braided tensor category of modular double representations of $U_q(sl(2,\mathbb{R}))$ is equivalent to a certain braided tensor category of unitary representations of the Virasoro algebra. Then the analogy with the compact case suggests that there should be the corresponding  braided tensor category for $\widehat{sl(2,\mathbb{R})}$ and it must be composed of an affine counterpart of the principal series representations.

A generalization of the classical analytic construction of the principal series representations to the affine Lie algebra 
$\widehat{sl(2,\mathbb{R})}$ presents substantial analytic difficulties still unresolved even in $\widehat{su(2)}$ case, in which there exists only a heuristic physical construction, known as WZW model. 
In mathematical literature there exists a construction of a unitary representation of $\widehat{sl(2,\mathbb{R})}$, see \cite{ggv}. However, this construction has a trivial central extension and does not seem to be relevant for the problem of constructing an equivalent tensor category with the ones in \cite{teschner}.

In this paper we combine an analytic and algebraic approaches to construction of representations of 
$\widehat{sl(2,\mathbb{R})}$, following the analogy with the quantum algebra case. The principal series representation of the modular double of $U_q(sl(2,\mathbb{R}))$ was realized by means of the algebra of the quantum plane, which is a quantum version of the 
$ax+b$-group \cite{ivan}. In our construction of $\widehat{sl(2,\mathbb{R})}$ representations we exploit the representations of the loop $ax+b$-group and the corresponding Lie algebra \cite{zeit}. The resulting formulas for generators at the first glance remind formulas for bosonization of the free field representation of highest weight modules \cite{gmm}. However, we encounter a problem of divergences in the correlation functions of the  $\widehat{sl(2,\mathbb{R})}$ generators that endangers 
the foundations of our construction: this is the major issue which makes our construction different from highest weight case. 
One of the central results of our paper is the method of eliminating divergencies so that the resulting "renormalized" correlators still satisfy the correct commutation relations of $\widehat{sl(2,\mathbb{R})}$ Lie algebra. Our approach was inspired by the explicit combinatorial formula for the correlators of affine Lie algebra 
in a highest weight representation discovered in \cite{zhu}. The formula is expressed in terms of Feynman-like diagrams, where the connected ones had the form of cycles.

In our case, one can write down the expression for regularized correlators in terms of Feynman-like diagrams, where the connected diagrams involve tree and 1-loop graphs. We propose, that there is an explicit way to renormalize 1-loop graphs, so that the commutation relations and the Hermicity condition within the correlator remain intact. So, our construction produces representations for $\widehat{sl(2,\mathbb{R})}$ with a bilinear Hermitian form, which depends on an infinite family of renormalization parameters. It is an important problem to find a set of parameters for which the resulting representation, determined by renormalized correlators is unitary, i.e. the pairing is nondegenerate and positive-definite. 

Following the analogy with the representation of the 
quantum group, one can generalize the construction of this article to the higher rank. In the quantum group case one had to use the representations of the multiple quantum planes \cite{gerasimov}, \cite{kashaev}, \cite{fip}, \cite{ivan2}, \cite{ivan3}. In our case it seems to be a similar technical problem: 
the role of quantum planes will be played by the loop $ax+b$-group. 

The structure of the article is as follows. In Section 2 we recall the structure of principal series representations 
of $sl(2,\mathbb{R})$, we reformulate it algebraically, indicating its relation to the representations of *-algebras $\mathcal{A}$,  $\mathcal{K}$, the two versions of $ax+b$-algebra. 
In Section 3 we discuss the unitary and nonunitary representations of affinized versions of  *-algebras $\mathcal{A}$,  $\mathcal{K}$, which we denote as $\widehat{\mathcal{A}}$,  $\widehat{\mathcal{K}}$ correspondingly.  
Section 4 is devoted to our main results, namely a construction of the $\widehat{sl(2,\mathbb{R})}$ representations. First we give formal expressions for the generators of $\widehat{sl(2,\mathbb{R})}$, built from the generators of $\widehat{\mathcal{A}}$,  $\widehat{\mathcal{K}}$ correspondingly, 
which satisfy the commutation relations of $\widehat{sl(2,\mathbb{R})}$ but which lead to divergent correlation functions for 
the $\widehat{\mathcal{A}}$,  $\widehat{\mathcal{K}}$-representations. We show how to obtain the well-defined correlation functions of the generators, which determine the representation. At first we show how it works for nonunitary representations of  $\widehat{\mathcal{A}}$, $\widehat{\mathcal{K}}$ as an example and then generalize the buildup to the unitary ones.\\

\noindent{\bf Acknowledgements.} We are grateful to H. Garland and A. Goncharov for illuminating discussions. A.M.Z. is indebted to the organizers of the Simons Workshop 2012, where this paper was partly written. I.B.F. is supported by 
NSF grant DMS-100163.

\section{A realization of continuous series for $sl(2,\mathbb{R})$}
\noindent {\bf 2.1. $\mathcal{A}$ and $\mathcal{K}$ algebras and their unitary representations.} 
In this section we present a construction of continuous series of $sl(2,\mathbb{R})$ in a special way, which will be convinient for generalizations to the loop case. We first construct certain versions of ax+b-algebras (affine transformations of a real line) which we call 
$\mathcal{A}$-algebra and $\mathcal{K}$-algebra: they have similar algebraic structure, but a different *-structure.
The $\mathcal{A}$-algebra is an algebra with three generators: $h,e^{\pm}$, so that  $ih,ie^{+}$ generate the Lie algebra of 
ax+b-group, however we make the b-subgroup generator $e^{+}$ to be invertible, so that the inverse element is $e^{-}$.
The commutation relations and the star structure are as follows:
\begin{eqnarray}
[h, e^{\pm}]=\pm ie^{\pm}, \quad e^{\pm}e^{\mp}=1, \quad h^*=h, \quad e^{{\pm}^*}=e^{\pm}.
\end{eqnarray}
The generators $h, \alpha^{\pm}$ of $\mathcal{K}$-algebra have similar commutation relations, but the star-structure is different:
 \begin{eqnarray}
 [h, \alpha^{\pm}]=\mp \alpha^{\pm}, \quad \alpha^{\pm}\alpha^{\mp}=1, \quad h^*=h, \quad {\alpha^{\pm}}^*=\alpha^{\mp}.
\end{eqnarray}
To construct the unitary representations of the above *-algebras, one can realize the generators $h,e^{\pm}$ as unbounded self-adjoint operators $i\frac{d}{dx}, e^{\pm x}$ in Hilbert space of square-integrable functions on a real line $L_2(\mathbb{R})$ as well as generators $h, \alpha^{\pm}$ as the operators $i\frac{d}{d\phi}, e^{\pm i\phi}$ in Hilbert space of square integrable functions on a circle $L_2(S^1)$. These generators acting on a specific vector 
in the appropriate Hilbert space generate a dense set. More explicitly one has the following.\\ 

\noindent {\bf Proposition 2.1.} {\it i) Let $D_A$ be the space spanned by the vectors $a\cdot v_0$, where $v_0=e^{-t x^2}\in L_2(\mathbb{R})$, $t>0$ and 
$a$ belongs to universal enveloping algebra of $\mathcal{A}$-algebra, such that the action of generators is realized as   
$h=i\frac{d}{dx}, e^{\pm}= e^{\pm x}$. Then $D_A$ is a dense set in $L_2(\mathbb{R})$.\\

\noindent ii) Let $D_K$ be the space spanned by the vectors $a\cdot v_0$, where $v_0=1\in L_2(S^1)$ and 
$a$ belongs to universal enveloping algebra of $\mathcal{K}$-algebra, such that the action of generators is realized as   
$h=i\frac{d}{d\phi}, \alpha^{\pm}= e^{\pm i \phi}$, where $\phi\in [0,2\pi]\cong S^1$ . Then $D_K$ is a dense set in $L_2(S^1)$.}\\

\noindent{\bf Proof.}  i) Let us consider the space spanned by vectors $h^nv_0$. This way one can obtain all functions of the form $P(x)e^{-tx^2}$, where $P(x)$ is any polynomial function. It is well-known that this set is dense in $L_2(\mathbb{R})$. Therefore  $D_A$ is dense in $L_2(\mathbb{R})$.\\
\noindent ii) In this case let us consider the space spanned by vectors $(\alpha^+)^n(\alpha^-)^m v_0$. This way we obtain the space of all trigonometric polynomials. It is dense in $L_2(S^1)$. Therefore $D_K$ is a dense set in $L_2(S^1)$.
\hfill $\blacksquare$.\\

\noindent {\bf 2.2. Continuous series of $sl(2,\mathbb{R})$ via $\mathcal{A}$ and $\mathcal{K}$ algebras.} 
 First of all, let us introduce some notations. We are interested in unitary representations of $sl(2,\mathbb{R})$ algebra, i.e. Lie algebra with generators $E,F, H$ such that 
 \begin{eqnarray}\label{sl2r}
&& [E,F]=H, \quad [H, E]=2E,\quad [H, F]=-2F\nonumber\\
&& E=-E^*, \quad F=-F^*, \quad H=-H^* .
 \end{eqnarray}
We remind that there are two standard relaizations of continuous series of $sl(2,\mathbb{R})$, one is related to inducing the representations from the diagonal subgroup of $sl(2,\mathbb{R})$, corresponding to the generator $H$, the other one is related to the inducing representations from maximally compact subgroup generated by $J^3=E-F$. It is convenient to introduce the following change of generators (which provides the correspondence between $sl(2,\mathbb{R})$ and $su(1,1)$):
\begin{eqnarray}
J^{\pm}=E+F\mp i H ,  
\end{eqnarray}
so that the relations \rf{sl2r} can be rewritten in the following way:
\begin{eqnarray}
&&[J^3, J^{\pm}]=\pm 2iJ^{\pm}, \quad 
[J^+, J^-]=-iJ^3, \nonumber \\
&&{J^+}^*=-J^-, \quad {J^3}^*=-J^3.
 \end{eqnarray}
Now we will give two realizations of $sl(2,\mathbb{R})$ via generators of $\mathcal{A}$- and $\mathcal{K}$-algebras.
Using $\mathcal{A}$-algebra one can write down the expressions for $E, F, H$:
\begin{eqnarray}\label{real1}
&&E=\frac{i}{2}(e^+h+he^+)+i\lambda e^+,\nonumber\\
&&F=-\frac{i}{2}(e^-h +he^-)+i\lambda e^-,\\
&&H=-2ih,\nonumber
 \end{eqnarray}
where $\lambda$ is a real parameter. 
Similarly, using $\mathcal{K}$-algebra, one can represent $J^3, J^{\pm}$ in a similar way:
\begin{eqnarray}\label{real2}
&&J^+=\frac{i}{2}(\alpha^+h+h\alpha^+)-\lambda \alpha^+\nonumber\\
&&J^-=\frac{i}{2}(\alpha^-h +h\alpha^-)+\lambda \alpha^-\\
&&J^3=2ih\nonumber
 \end{eqnarray}
 
Therefore the following theorem is valid.\\

\noindent {\bf Theorem 2.1.} {\it Let $\mathcal{D}_A$ be the space spanned by the vectors $a\cdot v_0$, where $v_0=e^{-t x^2}\in L_2(\mathbb{R})$ and 
$a$ belongs to universal enveloping algebra of $sl(2,\mathbb{R})$-algebra, such that the action of generators is realized as in \rf{real1}, so that  
$h=i\frac{d}{dx}, e^{\pm}= e^{\pm x}$. Then $\mathcal{D}_A$ is a dense set in $L_2(\mathbb{R})$ and it is a representation space for $sl(2,\mathbb{R})$.

\noindent ii) Let $\mathcal{D}_K$ be the space spanned by the vectors $a\cdot v_0$, where $v_0=1\in L_2(S^1)$ and 
$a$ belongs to universal enveloping algebra of $sl(2,\mathbb{R})$-algebra, such that the action of generators is realized as in \rf{real2}, so that $h=i\frac{d}{d\phi}, e^{\pm}= e^{\pm i \phi}$, where $\phi\in [0,2\pi]\cong S^1$. Then $\mathcal{D}_K$ is a dense set in $L_2(S^1)$ and it is representation space for $sl(2,\mathbb{R})$.} \\

\noindent{\bf Proof.}  i) Let us consider the space spanned by vectors $H^nv_0$. This way one can obtain all functions of the form $P(x)e^{-tx^2}$, where $P(x)$ is any polynomial function. It is well-known that this set is dense in $L_2(\mathbb{R})$. Therefore  $\mathcal{D}_A$ is dense in $L_2(\mathbb{R})$.

ii) In this case let us consider the space spanned by vectors $(J^+)^n(J^-)^m v_0$. This way we obtain the space of all trigonometric polynomials. It is dense in $L_2(S^1)$. Therefore $\mathcal{D}_K$ is a dense set in $L_2(S^1)$.
\hfill $\blacksquare$.\\

\section{Representations of $\widehat{\mathcal{A}}$ and $\widehat{\mathcal{K}}$}
\noindent {\bf 3.1. Definitions.} In this section we consider the following *-algebras, which we denote $\widehat{\mathcal{A}}$ and $\widehat{\mathcal{K}}$. The first algebra is close to the Lie algebra of a loop group associated with extended $\mathcal{A}$-group with properties, which is generated by $h_n$, $e^{\pm}_n$, $n\in \mathbb{Z}$,
so that the generating "currents" are:
\begin{eqnarray}
h(u)=\sum_{n\in\mathbb{Z}} h_{-n} e^{inu}, \quad e^{\pm}(u)=\sum_{n\in\mathbb{Z}} e^{\pm}_{-n}e^{inu} 
\end{eqnarray}
and obey the following commutation relations, expressed via generating functions:
\begin{eqnarray}\label{relcom}
&&[h(u), h(v)]=0, \quad [e^{\pm}(u), e^{\pm}(v)]=0, \quad [e^{\pm}(u), e^{\mp}(v)]=0,\nonumber\\
&&[h(u), e^{\pm}(v)]=\pm i\delta(u-v)e^{\pm}(v), \quad e^{\pm}(u)\cdot e^{\mp}(u)=1.
\end{eqnarray}
The *-structure is such that $h(u)$, $e^{\pm}(u)$ are Hermitian, i.e.
\begin{eqnarray}
h(u)^*=h(u),\quad e^{\pm}(u)^*=e^{\pm}(u)
\end{eqnarray}
The algebra $\widehat{\mathcal{K}}$ is generated by $h_n$, $\alpha^{\pm}_n$, $n\in \mathbb{Z}$,
it has similar commutation relations, but a different *-structure:
\begin{eqnarray}
&&h(u)=\sum_n h_{-n} e^{inu}, \quad \alpha^{\pm}(u)=\sum_n \alpha^{\pm}_{-n}e^{inu},\nonumber\\ 
&&[h(u), h(v)]=0, \quad [\alpha^{\pm}(u), \alpha^{\pm}(v)]=0, \quad [\alpha^{\pm}(u), \alpha^{\mp}(v)]=0,\nonumber\\
&&[h(u), \alpha^{\pm}(v)]=\mp \delta(u-v)\alpha^{\pm}(v), \quad \alpha^{\pm}(u) \alpha^{\mp}(u)=1\nonumber\\
&&h(u)^*=h(u),\quad \alpha^{\pm}(u)^*=\alpha^{\mp}(u).
\end{eqnarray}
After a short reminder of Gaussian integration, we will construct some representations of these algebras.\\

\noindent{\bf 3.2. Gaussian integration on Hilbert spaces.} 
In this subsection we recall a few basic facts and formulas, for more information, see e.g. \cite{daprata}, \cite{quo}. 
Suppose we have a real separable Hilbert space $H$ with the orthonormal basis $\{e^i\}$, $i\in \mathbb{N}$ and the pairing $\langle\cdot , \cdot\rangle$. Every element $x$ of this Hilbert space can be expressed as $x=\sum_i x_ie^i$. Let us introduce positive real numbers $\lambda_i$, $i\in \mathbb{N}$, so that 
 $\sum_i\lambda_i<\infty$. Then one can say that numbers $\lambda_i$ define a diagonal trace class operator $A$ on our Hilbert space. 
Than it appears possible to define the sigma-additive measure $d\mu_A$ on $H$ 
and heuristically express it as follows:
\begin{eqnarray}
d\mu_A(x)=(\sqrt{\det {2 \pi A}})^{-1}\cdot e^{-\frac{1}{2}\langle x,A^{-1}x \rangle}[dx],
\end{eqnarray}
 which can be thought as the infinite product of 1-dimensional Gaussian measures for each $i$: $d\mu_i=\sqrt{2\pi \lambda_i}^{-1}e^{-\lambda_i^{-1}x_i^2}$.

Since it is a sigma-additive measure, one can define the space of square-integrable functions $L_2(H, d\mu_A)$ with respect to it. One of the basic formulas
is the translational shift in the measure. Namely, if $b\in Im A$, then 
\begin{eqnarray}
\int f(x) d\mu_A(x)=\int f(x+b) e^{{-\frac{1}{2}\langle b,A^{-1}b \rangle}-\langle x,A^{-1}b \rangle}d\mu_A(x).
\end{eqnarray}
One can consider also an infinitesimal version of this formula. Making $b$ infinitesimal and parallel to $e_i$, we obtain 
that 
\begin{eqnarray}
\int D_i f(x)d\mu_A(x)=0, \quad D_i=\p_{x_i}-\lambda_i^{-1}x_i,
\end{eqnarray}
if $\p_{x_i} f(x)\in L_2(H, d\mu_A)$.  
It should be noted that the following monomials
\begin{eqnarray}\label{polexp}
(\prod^k_{i=1}\langle \alpha_i, x\rangle) e^{\langle \beta, x\rangle},
\end{eqnarray}
where $\alpha_i, \beta$ are the elements of complexified Hilbert space, are always integrable with respect to $d\mu_A$, moreover, they belong  to $L_2(H, d\mu_A)$. There is an explicit formula for the integral of the function \rf{polexp}, which can be derived from the simple result:
\begin{eqnarray}\label{gauss}
\int  e^{\langle \beta, x\rangle}d\mu_A(x)=e^{\frac{1}{2}\langle \beta, A\beta\rangle}.
\end{eqnarray}

\noindent{\bf 3.3. Construction of representations.} In order to construct representations of $\widehat{\mathcal{A}}$ and $\widehat{\mathcal{K}}$ one can use the Gaussian measure on a Hilbert space. 
Let us consider the Fourier series of a function from $L_2(S^1, \mathbb{R})$:
\begin{eqnarray}
x(u)=\sum_{n\in \mathbb{Z}}x_{-n}e^{inu}, \quad x_0\in \mathbb{R}, \quad x_n^*=x_{-n}.
\end{eqnarray}
Let us introduce two quadratic forms defining two types of trace-class operators on $L_2(S^1, \mathbb{R})$, which will determine the appropriate Gaussian measures:
\begin{eqnarray}\label{bil}
&& B_A(x,x)=\frac{1}{2}\sum_{n\ge 1}\xi_n^{-1}x_nx_{-n}+\xi_0^{-1}x_0^2,\nonumber\\
&& B_{K}(x,x)=\frac{1}{2}\sum_{n\ge 1}\xi_n^{-1}x_nx_{-n},
\end{eqnarray}
where $\xi_n>0$ for all $n$ and $\sum_n\xi_n<\infty$. 
The Gaussian measures we are interested in, heuristically can be expressed as follows:
\begin{eqnarray}
&& dw_{A}=(\sqrt{det(2\pi N_A)})^{-1}e^{-B_A(x,x)}dx_0\prod^{\infty}_{n=1}[\frac{i}{2}dx_n\wedge dx_{-n}],\nonumber \\
&& dw_{K}=(\sqrt{det(2 \pi N_{K})})^{-1}e^{-B_K(x,x)}d\phi\prod^{\infty}_{n=1}[\frac{i}{2}dx_n\wedge dx_{-n}],
\end{eqnarray}
where $N_A, N_{K}$  are trace-class diagonal operators determined by the quadratics forms \rf{bil}. Here as before the range for $\phi$ is $[0,2\pi]$. Literally the difference between two measures is that in the second one we compactified the zero mode $x_0$ on a circle with parameter $\phi$ for $dw_{K}$. The Hilbert spaces of square-integrable functions with respect to these measures are denoted in the following as $\mathcal{H}_A$ and $\mathcal{H}_K$. 

 Let us construct the unitary representations of $\widehat{\mathcal{A}}$ and $\widehat{\mathcal{K}}$-algebras in the Hilbert spaces $\mathcal{H}_A$ and $\mathcal{H}_K$, correspondingly. 
 
 We will explicitly define the operators, which will represent the generators. Let us start from  $\widehat{\mathcal{A}}$-algebra. First of all, let us extend the index of $\xi_n$ to all integers, so that $\xi_{-n}=\xi_n$, $n\in \mathbb{Z}$. We need the following differential operators: 
 \begin{eqnarray}\label{ab}
b_{-n}=i(\p_n-\xi_n^{-1}x_{-n}), 
\quad a_{-n}=i\p_n,
\end{eqnarray}
where $\p_{n}=\frac{\p}{\p x_n}$. These operators are formally conjugate to each other, 
\begin{eqnarray}
a_n^*=b_{-n},
\end{eqnarray}
considered on a certain dense set $D^l_A$, i.e. functions which are the sums of monomials 
\begin{eqnarray}
\prod^n_{k=1}\langle\mu_k,x\rangle\prod^m_{s=1} \langle \lambda_s, e^{\pm x}\rangle, 
\end{eqnarray}
where $\langle,\rangle$ is the standard $L_2(S^1, \mathbb{R})$ pairing and $\mu_k, \lambda_k$ are trigonometric polynomials. 
We define the operators $h_n$ as follows:
\begin{eqnarray}\label{ha}
h_{-n}=\frac{1}{2}(a_{-n}+b_{-n})=i(\p_n-\frac{1}{2}\xi_n^{-1}x_{-n}).
\end{eqnarray}
Hence, on a dense set $h_n^*=h_{-n}$, so that the current $h(u)=\sum_{n\in \mathbb{Z}}h_{-n}e^{inu}$
is Hermitian. 
Also, note that Hermitian currents $h(u), x(v)$ generate infinite-dimensional Heisenberg algebra:
\begin{eqnarray}
[h(u), x(v)]=i\delta(u-v).
\end{eqnarray}
We also define the currents $e^{\pm}(u)$: 
\begin{eqnarray}\label{ea}
e^{\pm}(u)=e^{\pm x(u)}.
\end{eqnarray}
It also follows that they satisfy the commutation relations \rf{relcom}. One can show that $e^{\pm x(u)}$ for any $u$ are Hermitian operators, considered on a dense set $D^l_A$. Therefore, this gives a unitary representation of $\widehat{\mathcal{A}}$-algebra. \\

\noindent {\bf Proposition 3.1.} {\it Let us consider the elements of the form $a\cdot v_0$, where $v_0\equiv 1\in \mathcal{H}_A$ and $a$ belongs to universal enveloping algebra of  $\widehat{\mathcal{A}}$ and the action of 
the generators is given by the formulas \rf{ea}, \rf{ha}. Then such elements generate a dense set 
$D^l_A$ in $\mathcal{H}_A$, which is a unitary representation of  $\widehat{\mathcal{A}}$-algebra.}\\

\noindent{\bf Proof.} Idea of the proof is similar to the one from Proposition 2.1. Let us consider a linear span of the following vectors: $h^{k_1}_{n_1}....h^{k_p}_{n_p}v_0$ for $k_i\in \mathbb{Z}_{\ge 0}$, $n_i\in \mathbb{Z}$. 
The resulting space will contain all polynomials in $x_n$, which is a well-known dense set (see e.g. \cite{daprata}) in $\mathcal{H}_A$. Therefore $D^l_A$ is dense in $\mathcal{H}_A$. \hfill $\blacksquare$

Similarly one can construct the unitary representations of $\widehat{\mathcal{K}}$-algebra. Let us consider the dense set $D^l_A$ in $\mathcal{H}_K$ of the following form:
\begin{eqnarray}
\prod^n_{k=1}\langle\mu_k,x\rangle\prod^m_{s=1} \langle \lambda_s, e^{\pm ix^c}\rangle, 
\end{eqnarray}
where $\lambda_i$ are trigonometric polynomials without the constant term and $x^c(u)=\phi+\sum_{n\neq 0}x_n e^{-inu}$.
The currents $h(u), \alpha^{\pm}(u)$ are defined by the following formulas:
\begin{eqnarray}\label{ag}
&& h(u)=\sum_{n\neq 0} i(\p_n-\frac{1}{2}\xi_n^{-1}x_{-n})e^{inu}+i\p_{\phi},\nonumber\\
&& \alpha^{\pm}(u)=e^{\pm ix^c(u)}.
\end{eqnarray}
It is possible to introduce operators $a_n, b_n$ and  define them by
the same formulas as in \rf{ab} for all $n\neq 0$. For $n=0$ we put $a_0=b_0=i\p_{\phi}$. 
This allows to formulate the following.\\

\noindent {\bf Proposition 3.2.} {\it Let us consider the elements of the form $a\cdot v_0$, where $v_0\equiv 1\in \mathcal{H}_K$ and $a$ belongs to universal enveloping algebra of  $\widehat{\mathcal{K}}$ and the action of 
the generators is given by the formulas \rf{ag}. Such elements generate a dense set 
$D^l_K$ in $\mathcal{H}_K$, which is a unitary representation of  $\widehat{\mathcal{K}}$-algebra.}\\

\noindent{\bf Proof.} The proof is ismilar to the one of Proposition 3.1. Let us consider a linear span of the following vectors: 
$h^{k_1}_{n_1}....h^{k_p}_{n_p}v_0$ for $k_i\in \mathbb{Z}_{\ge 0}$, $n_i\in \mathbb{Z}$. 
The resulting space will contain all polynomials of $x_n$, where $n\neq 0$. This would be a dense set in the Hilbert space with the Gaussian measure with out the $\phi$-variable. If one applies the action of $\alpha_0^{\pm}$ modes to elements of this space and considers a linear span then the resulting set is dense in the space of polynomials of 
$x_n$, $e^{\pm i\phi}$ where $n\neq 0$, 
which is a dense set (see e.g. \cite{daprata}) in $\mathcal{H}_K$. Therefore $D^l_K$ is dense in $\mathcal{H}_A$. \hfill $\blacksquare$ \\

\noindent {\bf 3.4. Correlators and normal ordering.} 
An important notion which is necessary for our construction is the correlator associated with the representation. 

 By the correlator of generators $T_1, ..., T_n$ of $\widehat{\mathcal{A}}$-algebra (resp. $\widehat{\mathcal{K}}$-algebra) 
 we mean the following expression:
 \begin{eqnarray}
 <T_1...T_n>\equiv \langle v_0, T_1...T_n v_0 \rangle,
 \end{eqnarray}
 where the pairing $\langle, \rangle $ is of the Hilbert space $\mathcal{H}_A$ (resp. $\mathcal{H}_K$), $v_0$ is the vector corresponding to the constant function $1$ and  $T_1, ..., T_n$ 
 are the generators $h_n, e^{\pm}_m$ (resp. $h_n, \alpha^{\pm}_m$). 
 
We remind that $h_n=\frac{1}{2}(a_n+b_n)$.  We note the following property which is the consequence of the properties of Gaussian integration.\\

\noindent {\bf Proposition 3.3.}
{\it The following correlators  
\begin{eqnarray}\label{van}
<T_1...T_na_k>, \quad 
<b_kT_1...T_n>
\end{eqnarray}
vanish for any generators $T_1, ...,T_n$.}
\\

\noindent{\bf Proof}. The correlator of the first type $<T_1...T_na_k>= \langle v_0, T_1...T_n a_kv_0 \rangle$ vanishes, because 
$a_k v_0=0$ for any $k$. The correlator  $<b_kT_1...T_n>$ vanishes because by complex conjugation it transforms in the correaltor of the first type.\hfill $\blacksquare$
\\

It is natural to call operators $a_n$  $annihilation$ operators and 
$b_n$  $creation$ operators. This allows us to define the normal ordering. Namely, when we write down the expression 
\begin{eqnarray}
:T_1...T_n:
\end{eqnarray}
for the product of $n$ generators, we reorder them in such a way that creation operators will be to the left and annihilation to the right. 

This procedure together with the vanishing of the correlators \rf{van} also gives an easy method to compute the correlators: by means of commutation relations of the generators, one can reduce the products $T_1, ...,T_n$ to the normally ordered expressions. Therefore, the result  will reduce to the correlators of the generators $e^{\pm}_n$ or $\alpha^{\pm}_n$. 
 We also note that in the case of $\widehat{\mathcal{K}}$-algebra the corresponding correlators are nonzero only if they have an equal number of $\alpha^+$ and $\alpha^-$ generators. 
\\
In order to compute the correlators of these generators it is easier to consider the appropriate currents instead of modes and use the Gaussian integration.\\

\noindent {\bf Proposition 3.4.}{\it One has the following expressions for correlation functions:    
\begin{eqnarray}\label{core}
&&\langle e_+(u_1)...e_+(u_n)e_-(v_1)...e_-(v_m)\rangle=\nonumber\\
&&\exp(\sum^n_{i<j; i, j=1}N_A(u_i,u_j)+\sum^m_{r<s; r, s=1}N_A(v_r,v_s)-\sum^n_{k=1}\sum^m_{l=1}N_A(u_k,v_l)+\nonumber\\
&&\frac{n+m}{2}N_A(0,0)),\\
&&\label{corealpha}\langle \alpha_+(u_1)...\alpha_+(u_n)\alpha_-(v_1)...\alpha_-(v_m)\rangle=\\
&&\delta_{n, m}\exp(-\sum^n_{i<j; i, j=1}N_K(u_i,u_j)-\sum^n_{i<j; i, j=1}N_K(v_i,v_j)+\sum^n_{k,l=1}N_K(u_k,v_l)\nonumber\\
&&+nN_K(0,0)),
\end{eqnarray}
where 
\begin{eqnarray}
&&N_A(u, v)=2\sum_{n\ge 0}cos(n(u-v))\xi_n, \nonumber\\ 
&&N_K(u, v)=2\sum_{n>0}cos(n(u-v))\xi_n.
\end{eqnarray}}

\noindent{\bf Proof.} The proof of this result follows from the formula \rf{gauss}. Let us show that in the case of $e^{\pm}$-generators, the case of $\alpha^{\pm}$ ones is similar. We constructed $e^{pm}$ generators so that  
$e_+(u_1)...e_+(u_n)e_-(v_1)...e_-(v_m)=\exp{(\sum_i\delta_{u_i}+\sum_j\delta_{v_j}, x)}$, where $\delta_w$ stands for the delta function at the point $w$. Then one can see that if one treats $\sum_i\delta_{u_i}+\sum_j\delta_{v_j}$ as 
$\beta$ from \rf{gauss} one obtains the formula \rf{core}. To show that this is not only formally true, one has to consider $\delta_w$ as a weak limit of $L_2$-functions, the so-called $\delta$-like sequence, use that $\delta$-like approximation for $e^{\pm}$-generators, take a gaussian integral according to \rf{gauss} and then take a limit. For further details see e.g. \cite{daprata}.\hfill $\blacksquare$.\\

We note here that in the case of correlators of the currents involving not only $e^{\pm}, \alpha^{\pm}$, but 
also $h(u)$, the resulting expression will consist of monomials \rf{corealpha} multiplied on certain product of delta-functions, coming from commutation relations of $\widehat{\mathcal{K}}$, $\widehat{\mathcal{A}}$-algebras.
\\

\noindent {\bf 3.5. Nonunitary representations.} In the previous section we have shown that using Gaussian measure, one can construct unitary representations of $\widehat{\mathcal{A}}$, 
$\widehat{\mathcal{K}}$-algebras.  However, it is possible to simplify those representations by making $a_k$ commute with $b_s$ for any $k$ and $s$. However, the resulting module will not be unitary, i.e. the pairing though nondegenerate will lose its positivity. 

Let us at first give the explicit description of such module for $\widehat{\mathcal{K}}$-algebra. For this we consider the vacuum 
vector $v_0$, such that $v_0$ is annihilated by $a_k$:
\begin{eqnarray}
a_kv_0=0.
\end{eqnarray}
Then the module which we will refer to as $\mathcal{V}_A$ is spanned by the following vectors:
\begin{eqnarray}
b_{m_1}...b_{m_s}e^{\pm}_{n_1}...e^{\pm}_{n_r}v_0,
\end{eqnarray}
 where $n_1, ..., n_r, m_1, ..., m_s\in \mathbb{Z}$. 
 Let us  begin to define the pairing with the postulation of
the points:
$e(u)=\sum_ne_ne^{inu}$ is a Hermitian current, $b(u)^*=a(u)$, so that $b(u)=\sum_n b_ne^{-inu}$, $a(u)=\sum_n a_ne^{-inu}$ and 
\begin{eqnarray}
[b(u), a(v)]=0.
\end{eqnarray}
The pairing is uniquely defined by the correlator of $e^{\pm}$ currents is given by \rf{core}. 
Similarly one can define the module $\mathcal{V}_K$ for $\widehat{\mathcal{K}}$ algebra
with the same conditions, just replacing $e^{\pm}$ with $\alpha^{\pm}$, as a result $\mathcal{V}_K$ is spanned by
\begin{eqnarray}
b_{m_1}...b_{m_s}\alpha^{\pm}_{n_1}...\alpha^{\pm}_{n_r}v_0.
\end{eqnarray}
One can define the pairing on $\mathcal{V}_K$ which is uniquely defined by the correlators of $\alpha^{\pm}$ currents \rf{corealpha}. 
Let us formulate this as a proposition.\\

\noindent {\bf Proposition 3.5.} {\it  The pairing defined above is Hermitean and nondegenerate. It gives a structure of nonunitary representations of the *-algebras $\widehat{\mathcal{A}}$, $\widehat{\mathcal{K}}$ on the spaces 
$\mathcal{V}_A$, $\mathcal{V}_K$  correspondingly.}\\

We note here that clearly, these representations are nonunitary, because you can easily find vectors $v\in \mathcal{V}_A, \mathcal{V}_K$, such that $\langle v,v\rangle=0$.

\section{Construction of representations for $\widehat{sl(2,\mathbb{R})}$}
\noindent {\bf 4.1. Regularization and the commutator.} In order to construct the $\widehat{sl(2,\mathbb{R})}$ representations via $\widehat{\mathcal{A}}$, $\widehat{\mathcal{K}}$, one has to find affine analogue of the formulas for $E,F, H$ and $J^3, J^{\pm}$. In the case of nontrivial central extension, the $\widehat{\mathcal{A}}$, $\widehat{\mathcal{K}}$ representations appear to be insufficient, one has to introduce the representation of the infinite-dimensional Heisenberg algebra, so that the generating current is $\rho(u)=\sum_{n\in \mathbb{Z}}\rho_ne^{-inu}$ and the commutation relations are:
\begin{eqnarray}
[\rho_n, \rho_m]=2\kappa n\delta_{n,-m},
\end{eqnarray}
where $\kappa\in \mathbb{R}_{>0}$. 
The irreducible module, so-called Fock module $F_{\kappa,p}$ of this algebra is defined as follows. We introduce a vector $vac_p$ with the property $\rho_n vac_p=0$, $p\in \mathbb{R}$, $n>0$ so that
\begin{eqnarray}
F_{\kappa,p}=\{\rho_{-n_1}...\rho_{-n_k}vac_p; \quad n_1,..., n_k>0,\quad  \rho_0 vac_p=p\cdot vac_p\}.
\end{eqnarray}
The Hermitian pairing is defined so that 
\begin{eqnarray}
\langle vac_p, vac_p \rangle=1, \quad \rho^*_n=\rho_{-n}
\end{eqnarray}
Another object required in this section is the regularized version of the 
$\rho$, $h$, 
$e^{\pm}$, 
$\alpha^{\pm}$ currents. Namely, for any $\varphi$, which stands for any of $\rho$, $h$, 
$x(u)$ or $x^c(u)$ we consider 
\begin{eqnarray}
\varphi(z,\bar z)=\sum_{n\ge 0}\varphi_n \bar z^n+\sum_{n> 0}\varphi_{-n} z^n,
\end{eqnarray}
where $z=re^{iu}$, so that $0<r\le 1$. We denote $e^{\pm}(z,\b z)\equiv e^{\pm x(z,\b z)}$ and $\alpha^{\pm}(z,\b z)\equiv e^{\pm i x^c(z,\b z)}$. 
The Wick theorem implies that the correlators of the regularized Heisenberg currents $\rho(z_1, \bar z_1), ...., \rho(z_n,\bar z_n)$ are finite as long as $0<|z_i|< 1$, i.e. 
the expressions
\begin{eqnarray}\label{heicor}
\langle vac_p, \rho(z_1,\bar z_1)...\rho(z_n,\bar z_n) vac_p\rangle
\end{eqnarray}
are finite. One can consider the limit  $|z_i|\to1$ in the sense of distributions, so that \rf{heicor} is the sum of products of distributions, i.e. 
\begin{eqnarray}\label{2pointf}
\langle vac_p, \rho(u_1)\rho(u_2)vac_p\rangle=\frac{2\kappa}{(1-e^{i(u_2-u_1-i0)})^2}+p^2.
\end{eqnarray}
Next we consider the following sets of composite regularized currents 
\begin{eqnarray}
&&J^{\pm}(z,\bar z)=\\
&&\frac{i}{2}(b(z,\bar z)\alpha^{\pm}(z,\bar z)+\alpha^{\pm}(z,\bar z)a(z,\bar z))\pm\kappa\p_u\alpha^{\pm}(z,\bar z)
\pm \rho(z,\bar z)\alpha^{\pm}(z,\bar z),\nonumber\\
&&J^3(z,\bar z)=2i h(z,\bar z)-2\kappa\alpha^-(z,\bar z)\p_u \alpha^+(z,\bar z),\nonumber
\end{eqnarray}
and 
\begin{eqnarray}
&&E(z,\bar z)=\\
&&\frac{i}{2}(b(z,\bar z)e^{+}(z,\bar z)+e^{+}(z,\bar z)a(z,\bar z))+i\kappa\p_ue^{+}(z,\bar z)
+ i\rho(z,\bar z)e^{+}(z,\bar z),\nonumber\\
&&F(z,\bar z)=\nonumber\\
&&-\frac{i}{2}(b(z,\bar z)e^{-}(z,\bar z)+e^{-}(z,\bar z)a(z,\bar z))+i\kappa\p_ue^{-}(z,\bar z)
+i \rho(z,\bar z)e^{-}(z,\bar z)\nonumber\\
&&H(z,\bar z)=-2i h(z,\bar z)+2i\kappa e^-(z,\bar z)\p_u e^+(z,\bar z).\nonumber
\end{eqnarray}
Moreover, we have the following Hermicity conditions:
\begin{eqnarray}
&&E(z,\bar z)^*=-E(z,\bar z), \quad F(z,\bar z)^*=-F(z,\bar z), \quad H(z,\bar z)^*=-H(z,\bar z),\nonumber\\
&&J^3(z,\bar z)^*=-J^3(z,\bar z), \quad J^{\pm}(z,\bar z)^*=-J^{\mp}(z,\bar z).
\end{eqnarray}
Then the following proposition is true.\\

\noindent {\bf Proposition 4.1.} {\it 
Let $\phi_k$ denote $E, F, H$ or $J^3, J^{\pm}$. Then the correlators
\begin{eqnarray}\label{gencor}
&&\langle \phi_{1}(z_1,\bar z_1)...\phi_{n}(z_n,\bar z_n)\rangle_p\equiv\nonumber\\
&&\langle v_0\otimes vac_p,\phi_{1}(z_1,\bar z_1)...\phi_{n}(z_n,\bar z_n)v_0\otimes vac_p\rangle,
\end{eqnarray}
 are well-defined for $0<|z_i|<1$. 
Moreover, if one of the currents is on the unit circle, i.e. considered at the point $|z|=1$, while all other are inside the unit circle,  correlator \rf{gencor} is also well-defined. 
}\\

\noindent {\bf Proof.} Consider the case when all parameters are inside the unit circle. 
Since $|z_1|, ..., |z_n|<1$ after the normal ordering procedure $\phi_{1}(z_1,\bar z_1)...\phi_{n}(z_n,\bar z_n)$ will be represented as normally ordered products with coefficients which are continuous functions of $u_1,..., u_n$, where 
$z_i=r_ie^{u_i}$. Therefore, the expression $\langle \phi_{1}(z_1,\bar z_1)...\phi_{n}(z_n,\bar z_n)\rangle_p$ will be given by the sum of correlators of exponentials $\alpha^{\pm}$ or $e^{\pm}$ which are well-defined even on the unit circle as we know from the formulas \rf{core}. In the case if only one of the currents $\phi_k(z_k, \bar z_k)$ is such that 
$|z_k|=e^{iu_k}$, the situation will not change. This is because during normal ordering procedure, the commutators of $a$- and $b$- parts of this 
generator with other ones again produce continuous functions of $u_1,..., u_n$, since all other $z_i$ are 
inside the unit circle and the normal ordering wouldn't produce any distributions\footnote{As we will see below, in the case when two or more currents $\phi_i$ have their arguments lying on the unit circle the correlator diverges. Our aim will be to eliminate these divergencies consistently with the algebraic structure.}.\hfill $\blacksquare$
\\

Let us make sense of the commutator of two currents on the unit circle as follows. We begin by considering the difference of two correlators 
\begin{eqnarray}
&&\langle \phi_{1}(z_1,\bar z_1)...\xi(w_1,\bar w_1)\eta(w_2,\bar w_2)...\phi_{n}(z_n,\bar z_n)\rangle_p-\nonumber\\
&&\langle \phi_{1}(z_1,\bar z_1)...\eta(w_2,\bar w_2)\xi(w_1,\bar w_1)...\phi_{n}(z_n,\bar z_n)\rangle_p.
\end{eqnarray} 
It is clearly well-defined, however, we want $w_1, w_2$ to lie on the unit circle. 
Then we have the following.\\

\noindent {\bf Proposition 4.2.} {\it 
The following limit
\begin{eqnarray}\label{comm}
&&\lim_{r_1,r_2\to 1}\big(\langle \phi_{1}(z_1,\bar z_1)...\xi(w_1,\bar w_1)\eta(w_2,\bar w_2)...\phi_{n}(z_n,\bar z_n)\rangle_p-\nonumber\\
&&\langle \phi_{1}(z_1,\bar z_1)...\eta(w_2,\bar w_2)\xi(w_1,\bar w_1)...\phi_{n}(z_n,\bar z_n)\rangle_p\big ),
\end{eqnarray}
exists in the sense of distributions, more specifically the answer will contain a linear combination of delta functions $\delta(u_1-u_2)$ and its derivatives. 
Here  $\phi_k$ stands for $E, F, H$ or $J^3, J^{\pm}$, $w_i=r_ie^{iu_i}$, $\xi(u_1)\equiv \xi(e^{iu_1},e^{-iu_1})$, $\eta(u_2)\equiv \eta(e^{iu_2},e^{-iu_2})$. }\\

\noindent {\bf Proof.}
The proof is based on the definition of the commutator \rf{comm} and the commuation relations of the 
Heisenberg algebra generated by $\rho(u)$ and $\widehat{\mathcal{A}}$, $\widehat{\mathcal{K}}$-algebras. 
We explicitly prove the $J^{\pm}, J^3$ commutation relations, for 
 $E,F, H$-currents can be obtained in a similar fashion. 
Let $z=re^{iu}$, $w=te^{iu}$. In the notations below we drop  the dependence on $\bar z$ variable in order to simplify the calculations. Let us start by computing the commutator of $J^+$ and $J^-$. Let us introduce the notation:
\begin{eqnarray}
j^{\pm}(z,\bar z)\equiv\frac{i}{2}(b(z,\bar z)\alpha^{\pm}(z,\bar z)+\alpha^{\pm}(z,\bar z)a(z,\bar z)).
\end{eqnarray}
Then 
\begin{eqnarray}
&& \lim_{r,t\to 1}[j^+(z,\bar z), j^-(\bar w, \bar w)]=-2h(u)\delta(u-v),\nonumber\\
&&\lim_{r,t\to 1}[j^+(z,\bar z), -\kappa\p_v\alpha^-(w,\bar w)]=-i\kappa\p_v(\delta(u-v)\alpha^-(v))\alpha^+(u),\nonumber\\
&&\lim_{r,t\to 1}[ \kappa\p_u\alpha^+(z,\bar z), j^-(w,\bar w)]=
i\kappa\p_u(\delta(u-v)\alpha^+(u))\alpha^-(v),\nonumber\\
&&-\lim_{r,t\to 1}[\rho(z)\alpha^+(z,\bar z), \rho(w)\alpha^-(w)]=2i\kappa\p_u\delta(u-v)\alpha^+(u)\alpha^-(v).
\end{eqnarray}
Summing all the terms and using the $\delta$-function properties we arrive to the desired commutation relation between $J^+$ and $J^-$ currents:
\begin{eqnarray}\label{jpm}
[J^+(u),J^-(v)]=i J^3(v)\delta(u-v)+4i\kappa\delta'(u-v).
\end{eqnarray}
Similarly,
\begin{eqnarray}
\lim_{r,t\to 1}[-2i h(z,\bar z),2\kappa\alpha^-(w,\bar w)\p_v \alpha^+(w,\bar w)]=-4i\kappa\p_u\delta(u-v).
\end{eqnarray}
This leads to commutation relations between $J^3$ currents, giving the proper central extension term:
\begin{eqnarray}\label{j3}
[J^3(u), J^3(v)]=-8i\kappa\delta'(u-v)
\end{eqnarray}
 Finally, the commutation relation formula for $J^3$ and $J^{\pm}$ follows from the formula
\begin{eqnarray}
\lim_{r,t\to 1}[2\kappa\alpha^-(z,\bar z)\p_u \alpha^+(z,\bar z), j^{\pm}(w, \bar w)]=\pm2i\kappa\p_v\delta(u-v)\alpha^{\pm}(v).
\end{eqnarray}
Therefore, 
\begin{eqnarray}\label{j3pm}
[J^3(u), J^{\pm}(v)]=\pm 2iJ^{\pm}(v)\delta(u-v).
\end{eqnarray}
\hfill$\blacksquare$\bigskip

We will denote the  expression \rf{comm} as follows:
\begin{eqnarray}
\langle \phi_{1}(z_1,\bar z_1)...[\xi(u_1),\eta(u_2)]...\phi_{n}(z_n,\bar z_n)\rangle_p. 
\end{eqnarray}
We discovered throught the proof of Proposition 4.2 (see \rf{jpm}, \rf{j3}, \rf{j3pm})  that the commutation relations between currents $J^3, J^{\pm}$ (similarly one can show that for $E,F, H$-currents), generate the 
$\widehat{sl(2,\mathbb{R})}$ algebra. However, it doesn't mean that correlators \rf{gencor} define the 
representation of $\widehat{sl(2,\mathbb{R})}$, since the currents in \rf{gencor} will be often not well defined if two or more of the arguments are considered on the unit circle. Let us summarize what do we have so far in the following theorem.
\\

\noindent {\bf Theorem 4.1.} {\it The commutator, defined by the formula \rf{comm}, of the currents $E, F, H$ or $J^3, J^{\pm}$ exists and satisfies  the commutation relations \footnote{Again, we draw attention that these formulas are valid and well-defined only under the correlator as in Proposition 4.2.} for 
$\widehat{sl(2,\mathbb{R})}$ algebra with the central charge $\kappa$:
\begin{eqnarray}
&&[E(u),F(v)]=H(v)\delta(u-v)-4i\kappa\delta'(u-v), \quad [H(u), H(v)]=8i\kappa\delta'(u-v),\nonumber\\
&& [H(u), E(v)]=2E(v)\delta(u-v),\quad [H(u), F(v)]=-2F(v)\delta(u-v)
\end{eqnarray}
and
\begin{eqnarray}
&&[J^+(u),J^-(v)]=i J^3(v)\delta(u-v)+4i\kappa\delta'(u-v), \nonumber\\
&& [J^3(u), J^3(v)]=-8i\kappa\delta'(u-v),\nonumber\\
&&[J^3(u), J^{\pm}(v)]=\pm 2iJ^{\pm}(v)\delta(u-v).
\end{eqnarray}}

As we mentioned above, 
even though we managed to define the commutator, based on the regularized commutators, this definition doesn't provide a representation, since the correlation functions of $E, F, H$ or $J^3, J^{\pm}$ do not exist, if more than one of the arguments lies on a circle. Then it is clear that the space from which we started, i.e. $\mathcal{H}_K\otimes F_p$ or  $\mathcal{H}_A\otimes F_p$ are not suitable to be spaces for $\widehat{sl(2,\mathbb{R})}$-module. However, we 
still have the regularized correlators which obey commutation relations. If we manage to eliminate divergencies in such a way that commutation relations and Hermicity conditions would be preserved, then the correlators will determine representation with the Hermitian bilinear form.  
In the next subsection we show a method, how to get rid of the  divergencies and redefine the correlator, so that it is well-defined when all the arguments are on a circle.\\

\noindent {\bf 4.2. Renormalization of correlators and construction of representation: nonunitary representations.} 
In this section we will renormalize the correlators with currents on a circle for the cases of $\mathcal{V}_A$ and $\mathcal{V}_K$ representations. In both cases the description is very similar, so we focus on $\mathcal{V}_K$ case and the correlators of 
$J^3, J^{\pm}$ currents. Generalization to $\mathcal{V}_A$ and $E,F, H$ goes along the same path. 

In order to renormalize, at first we have to understand what kind of divergencies we are dealing with. For that purpose it is convenient to write down the expression for correlator in the graphic form using 
Feynman-like diagrams. 

Recall, that the easiest way to compute the correlator for $\widehat{\mathcal{K}}$ is to reduce the whole expression to the normal ordered form, i.e. all the creation operators $b(z,\bar z)$ are moved to the left and all the annihilation operators 
$a(z,\bar z)$ to the right.

Once we move the creation operator $b(z,\bar z)$, which is a part of the certain generator at the point to the left (or the annihilation operator $a(z,\bar z)$ to the right) it may produce the following terms arising from the commutation with the 
$\alpha^{\pm}$-generator of $J^{\pm}$ currents on the right (on the left in the case of $a(z,\bar z)$) of the given generator:
\begin{eqnarray}\label{lines}
&&[a(z,\bar z), \alpha^{\pm}(w,\bar w)]=\mp \alpha^{\pm}(w,\bar w)\delta(z,w),\nonumber\\
&&[\alpha^{\pm}(w,\bar w), b(z,\bar z)]= \pm \alpha^{\pm}(z,\bar z)\delta(z,w),
\end{eqnarray}
where $\delta(z,w)=\sum_{n\ge 0} (z\bar w)^n+\sum_{n> 0} (\bar zw)^n$. If $z,w$ are on the circle, i.e. $z=e^{iu}, w=e^{iv}$, we get $\delta(z,w)=\delta(u-v)$.
We will depict every term of the form \rf{lines}, which we obtain during the normal ordering procedure, as a line from one vertex to another: 

\begin{eqnarray}
\begin{xy}
(20,0)*+{\bullet}="2";
(20,5)*+{}="0";
(30,5)*+{\delta(z,w)}="1";
(40,0)*+{\bullet}="3";
"2";"3"  **\dir{-}; ?(.65)*\dir{>};
\end{xy}
\end{eqnarray}
Here the initial vertex corresponds to the term containing creation/annihilation operator and the terminal vertex correspond to the term containing $\alpha^{\pm}$. The direction of the arrow (to the right or to the left) indicates whether it was annihilation or creation operator.

At the same time, each generator $J^{\pm}$ contributes terms like that
\begin{eqnarray}
\alpha^{\pm}(z,\bar z)a(z,\bar z) , \quad b(z,\bar z)\alpha^{\pm}(z,\bar z).
\end{eqnarray}
They enter the diagram as vertices with one outgoing line to the right or to the left and any amount of incoming lines: 
\begin{eqnarray}
\begin{xy}
(0,0)*+{}="1";
(20,0)*+{\bullet}="2";
(20,3)*+{-}="-";
(20,5)*+{+}="+";
(40,0)*+{}="3";
(3,6)*+{}="1t";
(10,12)*+{}="1tt";
(3,-6)*+{}="1b";
(10,-12)*+{}="1bb";
 "1";"2" **\dir{-}; ?(.65)*\dir{>};
{\ar@{->} "2";"3"};
"1t";"2"  **\dir{-}; ?(.65)*\dir{>};
"1tt";"2"  **\dir{-}; ?(.65)*\dir{>};
"1b";"2"  **\dir{-}; ?(.65)*\dir{>};
"1bb";"2"  **\dir{-}; ?(.65)*\dir{>};
(60,0)*+{}="4";
(80,0)*+{\bullet}="5";
(80,3)*+{-}="-";
(80,5)*+{+}="+";
(100,0)*+{}="6";
(97,6)*+{}="6t";
(90,12)*+{}="6tt";
(97,-6)*+{}="6b";
(90,-12)*+{}="6bb";
 "5";"6" **\dir{-}; ?(.35)*\dir{<};
{\ar@{->} "5";"4"};
"5";"6t" **\dir{-}; ?(.35)*\dir{<};
"5";"6tt"  **\dir{-}; ?(.35)*\dir{<};
"5";"6b"  **\dir{-}; ?(.35)*\dir{<};
"5";"6bb"  **\dir{-}; ?(.35)*\dir{<};
\end{xy}
\end{eqnarray}
Signs $\pm$ over the vertices correspond to $\alpha^{\pm}$, while we neglect for simplicity the dependence on $z, \bar z$ variables. 
The incoming line in the vertex forms when during the normal ordering procedure the creation/annihilation operators from other vertices leave commutator term \rf{lines} with $\alpha^{\pm}$ at a given vertex. The outgoing line forms when the creation/annihilation operator of a given vertex leaves a commutator term $\alpha^{\pm}$ from another vertex. 

On the other hand, $J^{\pm}$ also have terms of the form 
\begin{eqnarray}\label{term}
\kappa\p_u\alpha^{\pm}(z,\bar z), \quad \rho(z,\bar z)\alpha^{\pm}(z,\bar z)
\end{eqnarray}
According to the strategy formulated, we denote the first term from \rf{term} as a terminal vertex for graphs, because it contains only $\alpha^{\pm}$: 

\begin{eqnarray}
\begin{xy}
(60,0)*+{}="4";
(80,0)*+{\bullet}="5";
(80,3)*+{-}="-";
(80,5)*+{+}="+";
(100,0)*+{}="6";
(97,6)*+{}="6t";
(90,12)*+{}="6tt";
(97,-6)*+{}="6b";
(90,-12)*+{}="6bb";
 "5";"6" **\dir{-}; ?(.35)*\dir{<};
"5";"6t" **\dir{-}; ?(.35)*\dir{<};
"5";"6tt"  **\dir{-}; ?(.35)*\dir{<};
"5";"6b"  **\dir{-}; ?(.35)*\dir{<};
"5";"6bb"  **\dir{-}; ?(.35)*\dir{<};
\end{xy}
\end{eqnarray}

However the second term from \rf{term} is composite: it has contribution from  
$\alpha^{\pm}(z,\bar z)$ and $\rho(z,\bar z)$. We consider this term as a terminal vertex for the graphs coming from the normal ordering of $a,b$-operators, however, we should add one outgoing "wavy" line, corresponding to the normal ordering in the Fock space: 
 
 \begin{eqnarray}
\begin{xy}
(60,0)*+{}="4";
(80,0)*+{\bullet}="5";
(80,3)*+{-}="-";
(80,5)*+{+}="+";
(100,0)*+{}="6";
(97,6)*+{}="6t";
(90,12)*+{}="6tt";
(97,-6)*+{}="6b";
(90,-12)*+{}="6bb";
{\ar@{~} "5";"4"};
 "5";"6" **\dir{-}; ?(.35)*\dir{<};
"5";"6t" **\dir{-}; ?(.35)*\dir{<};
"5";"6tt"  **\dir{-}; ?(.35)*\dir{<};
"5";"6b"  **\dir{-}; ?(.35)*\dir{<};
"5";"6bb"  **\dir{-}; ?(.35)*\dir{<};
\end{xy}
\end{eqnarray}
 
Each wavy line produces the term from 2-point correlator \rf{2pointf} and may appear only once in the connected graph: it is so, because of the combinatorial formula, expressing the Fock space n-point correlator via 2-point correlators.

Finally, there are terms coming from $J^3$-generator, namely
\begin{eqnarray}
&&\label{ab2}\frac{1}{2}a(z,\bar z), \quad  \frac{1}{2}b(z,\bar z),\\  
&&\label{terminalj3}2\kappa\alpha^-(z,\bar z)\p_u \alpha^+(z,\bar z)
\end{eqnarray}
According to our conventions, we denote first two elements \rf{ab2} as initial vertices with one outgoing line to the right and to the left correspondingly:
\begin{eqnarray}
\begin{xy}
(20,0)*+{\bullet}="2";
(20,5)*+{0}="0";
(40,0)*+{}="3";
"2";"3"  **\dir{-}; ?(.65)*\dir{>};
(60,0)*+{}="4";
(80,0)*+{\bullet}="5";
(80,5)*+{0}="0";
 "4";"5" **\dir{-}; ?(.35)*\dir{<};
\end{xy}
\end{eqnarray}
Zeros over the vertices represent the absence of $\alpha^{\pm}$ contribution from these terms.  
 Finally, \rf{terminalj3} will correspond to the terminal vertex, because it has neither creation or annihilation operators:  
\begin{eqnarray}
\begin{xy}
(60,0)*+{}="4";
(80,0)*+{\bullet}="5";
(80,5)*+{0}="0";
(100,0)*+{}="6";
(97,6)*+{}="6t";
(90,12)*+{}="6tt";
(97,-6)*+{}="6b";
(90,-12)*+{}="6bb";
 "5";"6" **\dir{-}; ?(.35)*\dir{<};
"5";"6t" **\dir{-}; ?(.35)*\dir{<};
"5";"6tt"  **\dir{-}; ?(.35)*\dir{<};
"5";"6b"  **\dir{-}; ?(.35)*\dir{<};
"5";"6bb"  **\dir{-}; ?(.35)*\dir{<};
\end{xy}
\end{eqnarray}
However, in this exceptional case we will consider every line as a sum of two terms contributing $\delta(z,w)$, because there are two $\alpha$'s in this vertex corresponding to the term $\alpha^-(z,\bar z)\p_u\alpha^+(z,\bar z)$, so we consider contribution of the commutator with each of them. 
Moreover, formally on the circle (if such a limit exists), there can be only one incoming line into this vertex, 
 because of the commutation relations:
\begin{eqnarray}
[a(u), \alpha^+(v)\p_v \alpha^-(v)]=[b(u), \alpha^+(v)\p_v \alpha^-(v)]=
-\delta'(u-v).
\end{eqnarray}
Here are two sample diagrams 

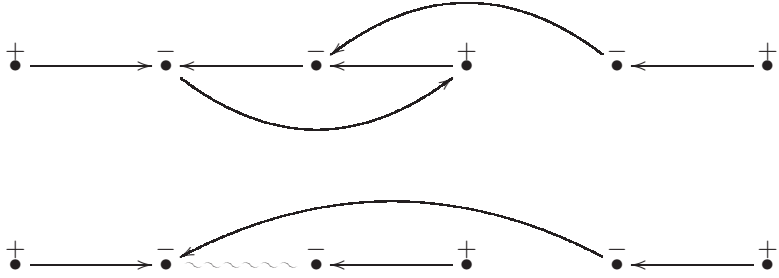
\begin{eqnarray}
\begin{xy}
(0,12)*+{+}="A+";
(20,12)*+{-}="B-";
(40,12)*+{-}="C-";
(60,12)*+{+}="D+";
(80,12)*+{-}="E-";
(0,10)*+{\bullet}="A";
(20,10)*+{\bullet}="B";
{\ar@{->} "A";"B"};
(40,10)*+{\bullet}="C";
(60,10)*+{\bullet}="D";
{\ar@{->} "C";"B"};
{\ar@/_2pc/ "B";"D"};
{\ar@{->} "D";"C"};
(80,10)*+{\bullet}="E";
{\ar@/_2pc/ "E";"C"};
(100,10)*+{\bullet}="F";
(100,12)*+{+}="F+";
{\ar@{->} "F";"E"}
\end{xy}\nonumber
\end{eqnarray}

\begin{eqnarray}
\begin{xy}
(0,12)*+{+}="A+";
(20,12)*+{-}="B-";
(40,12)*+{-}="C-";
(60,12)*+{+}="D+";
(80,12)*+{-}="E-";
(0,10)*+{\bullet}="A";
(20,10)*+{\bullet}="B";
{\ar@{->} "A";"B"};
(40,10)*+{\bullet}="C";
(60,10)*+{\bullet}="D";
{\ar@{~} "C";"B"};
{\ar@{->} "D";"C"};
(80,10)*+{\bullet}="E";
(100,10)*+{\bullet}="F";
(100,12)*+{+}="F+";
{\ar@/_2pc/ "E";"B"}
{\ar@{->} "F";"E"}
\end{xy}\nonumber
\end{eqnarray}

which belong to the expansion of the 6-point correlator
\begin{eqnarray}
\langle J^+(z_1,\bar z_1)J^-(z_2,\bar z_2)J^-(z_3,\bar z_3)J^+(z_4,\bar z_4)J^-(z_5,\bar z_5)J^+(z_6,\bar z_6)\rangle
\end{eqnarray}

After it was described how all terms from the generators fit into graphs describing the correlator, we are ready to formulate the first important statement about such graphical representation. \\

\noindent {\bf Proposition 4.3.} {\it Every connected graph contributing to the correlator 
\begin{eqnarray}\label{corgen}
\langle \phi_1(z_1,\bar z_1)....\phi_n(z_n,\bar z_n)\rangle_p,
\end{eqnarray}
where $\phi_i=J^3, J^{\pm}$, has at most one loop.}\\

\noindent {\bf Proof.} In order to see that, one has to analyze the structure of the vertices. 
Namely, each vertex have possibly many incoming lines, but at most one outgoing line. Let us consider the graph with one loop. It means that all outgoing lines of the vertices in this loop are included in the loop and all other lines are incoming.  Therefore none of these vertices can participate in a different loop, since all outgoing lines are included in the first loop. Let us assume that there is another loop with different set of vertices in the same graph. However, since for all vertices participating all outgoing lines are "circulating" inside the loop there is no possibility to make a connected graph out of two loops. 
\hfill$\blacksquare$\bigskip
\\
However, we see that each loop diagram, though  producing a well-defined expression for currents inside the circle, 
in the limit $|z_i|\to 1$ produce a divergence of the type $\delta(0)$, since every line in the loop produces delta-function. 
It is easy to see, say in the case of correlator 
\begin{eqnarray}
\langle J^+(z_1,\bar z_1)J^-(z_2,\bar z_2)\rangle_p
\end{eqnarray}
When the creation operator from $J^-$ contributes the commutator term \rf{lines} with the exponent in $J^+$ and at the same time, annihilation operator from 
$J^-$ contributes the commutator term \rf{lines} with the exponent in $J^+$, the following diagram is produced: 
\begin{eqnarray}
\xymatrix{
\bullet \ar@/_1pc/[r] &
\bullet   \ar@/_1pc/[l]}
\end{eqnarray}
\\
\noindent and leads to the divergent term $\delta(u_1-u_2)\cdot\delta(u_1-u_2)$ when considered on the circle. 

The simplest way to regularize correlators to preserve commutation relations is to throw away all the loop diagrams.  This allows to make sense of the correlators on a circle. \\

\noindent {\bf Theorem 4.2.a.} {\it Let us consider the renormalized correlators
 \begin{eqnarray}\label{rencor}
\langle \phi_1(z_1,\bar z_1)....\phi_n(z_n,\bar z_n)\rangle_p^R,
\end{eqnarray}
where $\phi_i=J^3, J^{\pm}$, which contain only the tree graph contributions to \rf{corgen}.
Then the limit ${r_k\to 1}$ (where $z_k=r_ke^{iu_k}$) of \rf{rencor} exists as a distribution. Moreover, the commutation relation between $J^3, J^{\pm}$ under the renormalized correlator reproduces $\widehat{sl(2,\mathbb{R})}$. 
Therefore, correlators \rf{rencor} considered on a circle define a module for $\widehat{sl(2,\mathbb{R})}$ with a Hermitian bilinear form.}\\

\noindent {\bf Proof.} 
To prove this theorem, it is enough to look at relevant contibutions for the commutators. It appears that they all come from initial/terminal vertices and their closest neighbors. For example, let us consider diagram of this sort: 
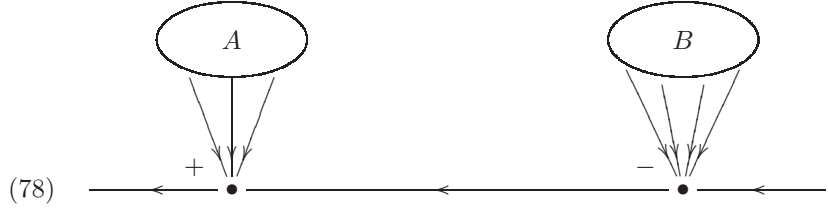
\begin{eqnarray}\label{diag1}
\begin{xy}
(0,0)*+{}="a";
(20,0)*+{\bullet}="b";
(80,0)*+{\bullet}="c";
(100,0)*+{}="d";
(15,3)*+{+}="+";
(75,3)*+{-}="-";
"a";"b"  **\dir{-}; ?(.45)*\dir{<};
"b";"c"  **\dir{-}; ?(.45)*\dir{<};
"c";"d"  **\dir{-}; ?(.45)*\dir{<};
(14,16)*+{}="b1";
(20,16)*+{}="b2";
(26,16)*+{}="b3";
"b1";"b"  **\dir{-}; ?(.75)*\dir{>};
"b2";"b"  **\dir{-}; ?(.75)*\dir{>};
"b3";"b"  **\dir{-}; ?(.75)*\dir{>};
(20,20)*+{A}="A";
(20,20)*\xycircle(10,5){};
(70,20)="e1";
(80,20)*+{B}="B";
(90,20)="e3";
(80,20)*\xycircle(10,5){};
(72,17)*+{}="c1";
(77,15)*+{}="c2";
(83,15)*+{}="c3";
(88,17)*+{}="c4";
"c1";"c"  **\dir{-}; ?(.75)*\dir{>};
"c2";"c"  **\dir{-}; ?(.72)*\dir{>};
"c3";"c"  **\dir{-}; ?(.72)*\dir{>};
"c4";"c"  **\dir{-}; ?(.75)*\dir{>};
\end{xy}
\end{eqnarray} 
\\
corresponding to the correlator 
\begin{eqnarray}\label{diag2}
\langle\dots J^+(z_1,\bar z_1)J^-(z_2,\bar z_2)\dots\rangle_p.
\end{eqnarray} 
There is a diagram with equal contribution 
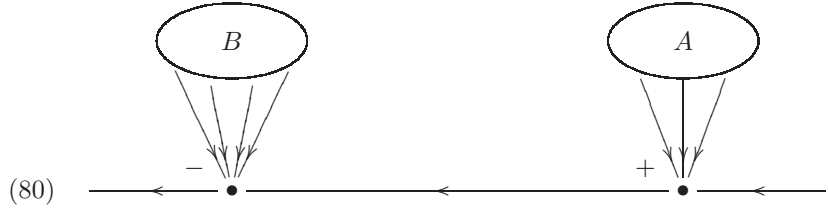
\begin{eqnarray}\label{diag2'}
\begin{xy}
(0,0)*+{}="a";
(20,0)*+{\bullet}="b";
(80,0)*+{\bullet}="c";
(100,0)*+{}="d";
(15,3)*+{-}="-";
(75,3)*+{+}="+";
"a";"b"  **\dir{-}; ?(.45)*\dir{<};
"b";"c"  **\dir{-}; ?(.45)*\dir{<};
"c";"d"  **\dir{-}; ?(.45)*\dir{<};
(12,17)*+{}="b1";
(17,15)*+{}="b2";
(23,15)*+{}="b3";
(28,17)*+{}="b4";
"b1";"b"  **\dir{-}; ?(.75)*\dir{>};
"b2";"b"  **\dir{-}; ?(.72)*\dir{>};
"b3";"b"  **\dir{-}; ?(.72)*\dir{>};
"b4";"b"  **\dir{-}; ?(.75)*\dir{>};
(20,20)*+{B}="B";
(20,20)*\xycircle(10,5){};
(70,20)="e1";
(80,20)*+{A}="A";
(90,20)="e3";
(80,20)*\xycircle(10,5){};
(74,16)*+{}="c1";
(80,16)*+{}="c2";
(86,16)*+{}="c3";
"c1";"c"  **\dir{-}; ?(.75)*\dir{>};
"c2";"c"  **\dir{-}; ?(.75)*\dir{>};
"c3";"c"  **\dir{-}; ?(.75)*\dir{>}
\end{xy}
\end{eqnarray}
\\
from the correlator
\begin{eqnarray}
\langle\dots J^-(z_2,\bar z_2)J^+(z_1,\bar z_1)\dots\rangle_p.
\end{eqnarray}
Therefore, they cancel each other when we take a commutator of $J^+(u_1), J^-(u_2)$, since the middle line produces $\delta(u_1-u_2)$. One can understand it in the following way: the commutator 
$[J^+(u_1), J^-(u_2)]$ involves $J^3$-term and 
$\delta'$-term, 
$J^3$-term contains only initial and terminal vertices, therefore the diagrams, which contribute to this commutator should have also $J^+$ or $J^-$ as initial or terminal vertices. As for the commutators 
of $J^3$ with $J^{\pm}$, the picture above indicates that only terminal/initial vertices and their closest neighbors contribute, because again, all terms in $J^3$ are depicted as initial/terminal vertices. Therefore, it is clear that one can cancel all loop contributions to the correlators and still the commutation relations of Theorem 4.1.will be valid. However, if all the loops are eliminated, one can consider the limit $|z_i|\to 1$, putting all currents on a circle, since all divergent graphs are gone. Therefore, we have the well-defined correlators of the generators of $\widehat{sl(2,\mathbb{R})}$ Lie algebra. One can see  
that the these correlators define a Hermitian bilinear form, because after the conjugation one can obtain one-to-one correspondence between tree graphs in conjugated correlators.
\hfill$\blacksquare$\bigskip
\\

It appears that the theorem above can be generalized: instead of eliminating loops completely, one can renormalize them, i.e. eliminate the divergence of the type $\delta(0)$. Namely, one can associate a real number $\mu_k$ with every loop which is the contribution of $k$ $J^{\pm}$-currents at the points $z_1,..., z_k$. This number $\mu_k$ enters the loop in the following way. 
Suppose currents $J^{\pm}(z_i)$ appear in the correlator in the given order from left to right. Our loop contains the following term
\begin{eqnarray}\label{delta1}
\delta (z_1,z_2)\cdot\delta(z_2,z_3)\dots\delta(z_{k-1}, z_k)\cdot \delta(z_{k}, z_1),
\end{eqnarray}
where $\delta (z,w)=\sum_{n\ge 0}(z\bar w)^n+\sum_{n< 0}(w\bar z)^n$. In the limit $|z_k|\to 1$ we will substitute it with  
the expression 
\begin{eqnarray}\label{delta2}
\mu_k\cdot \delta(u_1-u_2)\cdot\delta(u_2-u_3)\dots\delta(u_{k-1}-u_k),
\end{eqnarray}
where we remind that $z_k=r_ke^{iu_k}$. 
This procedure again does affect neither commutation relations inside the correlator nor Hermicity condition for the resulting bilinear form, because loops do not participate in commutation relations. Therefore,  we have a new version of the theorem above.
\\

\noindent {\bf Theorem 4.2.b.} {\it Let us consider renormalized correlators
 \begin{eqnarray}\label{rencormu}
\langle \phi_1(z_1,\bar z_1)....\phi_n(z_n,\bar z_n)\rangle_p^{R, \{\mu_n\}},
\end{eqnarray}
where $\phi_i=J^3, J^{\pm}$, such that the loop contributions to \rf{corgen} are replaced by their renormalized analogues with the family of arbitrary real parameters $\{\mu_n\}$. 
Then the limit ${r_k\to 1}$ (where $z_k=r_ke^{iu_k}$) of \rf{rencormu} exists and the commutation relations between $J^3, J^{\pm}$ under the renormalized correlator reproduce Lie algebra $\widehat{sl(2,\mathbb{R})}$, and correlators \rf{rencormu} considered on a circle define a module for $\widehat{sl(2,\mathbb{R})}$ with a Hermitian bilinear form.}\\

\noindent {\bf Proof.} The proof goes along the same lines as in the part $a$ of the theorem. When we consider the loop contribution to the commutator, there will be two loop diagrams, each emerging from one of the terms precisely as in the pictures \rf{diag1}, \rf{diag2} which again will cancel each other, because of the delta function line between two vertices and because there always is one incoming line and outgoing line for each vertex associated with the currents, participating in the commutator in these two terms. Therefore, the commutation relations are unaffected in the presence of loops. It is easy to see that the same applies to the Hermicity condition.
\hfill$\blacksquare$\bigskip
\\

\noindent{\bf 4.3. Renormalization of correlators and construction of representations: unitary representations.}
In the previous section, we have shown how to construct a module with a Hermitian pairing for $\widehat{sl(2,\mathbb{R})}$, starting from the nonunitary representations $\mathcal{V}_K$ of $\widehat{\mathcal{K}}$. In this section, we will do the same with the use of unitary module $D^l_K$, which is a dense subset in $H_K$. Needless to say, similar procedure will apply in the case of $\widehat{\mathcal{A}}$ and 
${D^l_A}$.
There is a major difference between the nonunitary and unitary case. The creation and annihilation operators $a_n$ and $b_n$ do not commute in this case, namely
\begin{eqnarray}
[a(z,\bar z), b(w,\bar w)]=D(z, w), \quad D(z,w)=\sum_{n>0}\xi_n^{-1}(z^n\bar w^n+\bar z^n w^n).
\end{eqnarray}
On a circle it is not a well-defined distribution, because of the behavior of $\xi_n$; however, it is well-defined on trigonometric polynomials and this is what we need to define the correlator of the modes of $J^{\pm}$.  Inside the circle, however, we 
can choose regularization parameters $z,w$ in such a way that $D(z,w)$ converges.   
When we compute the correlator
\begin{eqnarray}\label{corgen2}
\langle \phi_1(z_1,\bar z_1)....\phi_n(z_n,\bar z_n)\rangle_p,
\end{eqnarray}
which is now based on tensor product of Fock module and $D^l_K$, 
if $D(z,w)$ appears after we interchange $a, b$ operators, we indicate that with a dotted line. In other words, all the terms which contributed to vertices with solid outgoing lines in the case of nonunitary module, contribute also as vertices with outgoing dotted line and incoming solid lines. So, this is a set of new vertices for $J^{\pm}$ currents:
\begin{eqnarray}\label{dotted}
\begin{xy}
(0,0)*+{}="1";
(20,0)*+{\bullet}="2";
(20,3)*+{-}="-";
(20,5)*+{+}="+";
(40,0)*+{}="3";
(3,6)*+{}="1t";
(10,12)*+{}="1tt";
(3,-6)*+{}="1b";
(10,-12)*+{}="1bb";
 "1";"2" **\dir{-}; ?(.65)*\dir{>};
{\ar@{--} "2";"3"};
"1t";"2"  **\dir{-}; ?(.65)*\dir{>};
"1tt";"2"  **\dir{-}; ?(.65)*\dir{>};
"1b";"2"  **\dir{-}; ?(.65)*\dir{>};
"1bb";"2"  **\dir{-}; ?(.65)*\dir{>};
(60,0)*+{}="4";
(80,0)*+{\bullet}="5";
(80,3)*+{-}="-";
(80,5)*+{+}="+";
(100,0)*+{}="6";
(97,6)*+{}="6t";
(90,12)*+{}="6tt";
(97,-6)*+{}="6b";
(90,-12)*+{}="6bb";
 "5";"6" **\dir{-}; ?(.35)*\dir{<};
{\ar@{--} "5";"4"};
"5";"6t" **\dir{-}; ?(.35)*\dir{<};
"5";"6tt"  **\dir{-}; ?(.35)*\dir{<};
"5";"6b"  **\dir{-}; ?(.35)*\dir{<};
"5";"6bb"  **\dir{-}; ?(.35)*\dir{<};
\end{xy}
\end{eqnarray}
and $J^3$ current: 
\begin{eqnarray}
\begin{xy}
(20,0)*+{\bullet}="2";
(20,5)*+{0}="0";
(40,0)*+{}="3";
"2";"3"  **@{--}; 
(60,0)*+{}="4";
(80,0)*+{\bullet}="5";
(80,5)*+{0}="0";
 "4";"5" **@{--}; 
\end{xy}
\end{eqnarray}
One can generalize the results of the Theorem 4.2 in the case of the unitary modules. The renormalization procedure remains the same with arbitrary parameters $\mu_n$, except that in this case we will have some extra tree graphs  containing dotted lines, which may lead to divergencies, which we have to eliminate. \\

\noindent {\bf Theorem 4.3.} {\it Let us consider renormalized correlators
 \begin{eqnarray}\label{rencor2}
\langle \phi_1(z_1,\bar z_1)....\phi_n(z_n,\bar z_n)\rangle_p^{R, \{\mu_n\}}.
\end{eqnarray}
Here $\phi_i=J^3, J^{\pm}$. The loop contributions to \rf{corgen2} are replaced by their renormalized analogues with a  family of arbitrary real parameters $\{\mu_n\}$. The tree graphs, containing dotted lines are eliminated unless there is equal number of plus and minus vertices on one side of the dotted line \rf{dotted}. 
Then the limit ${r_k\to 1}$ (where $z_k=r_ke^{iu_k}$) of \rf{rencor} exists and the commutation relations between $J^3, J^{\pm}$ under the renormalized correlator reproduce Lie algebra $\widehat{sl(2,\mathbb{R})}$. The correlators \rf{rencor2} considered on a circle define the module for $\widehat{sl(2,\mathbb{R})}$ with the Hermitian bilinear form.}\\

\noindent {\bf Proof. }First, we make the following useful observation: there can be only one dotted line in the connected graph, and there is no loops in that graph. Alas, in most of cases these diagrams diverge on a circle because of the correlator between 
$\alpha^+,\alpha^-$, which, combined with $D(e^{iu},e^{iv})$ leads to divergence due to the product with $N_K(u,v)$. The basic example of such sort is the 2-point correlator between $J^+$ and $J^-$. Therefore, we will eliminate all of such graphs containing dotted lines so that the commutation relations remain intact. 

However, there is still a small amount of these diagrams, which leads to finite expressions and which is important for keeping the commutation relations. In order to describe them, let us separate vertices into two groups: those connected through solid path with first vertex with outgoing dotted line and those connected with another one. Note that if the number of 
pluses and minuses on vertices in any of these two groups is equal, then the expression corresponding to such a graph is well defined. Actually, every solid line gives a $\delta$-function, and since there is an equal number of $\alpha^+, \alpha^-$ on one side, then there will be no divergence arising from $N_K(u,v)$ and $D(e^{iu},e^{iv})$ as described above. All these graphs contribute to commutation relation between $J^+$ and $J^-$, which leads to emergence of $J^3$ current and the resulting graph with terminal/initial dotted line, involving $J^3$. 

In order to renormalize graphs with loops, one can follow the same route as in the previous subsection. 
In other words, one can renormalize all the loop diagrams in a similar way, using the formulas \rf{delta1}, \rf{delta2}. \hfill $\blacksquare$.


\begin{thebibliography}{25}
\bibitem{daprata}G. Da Prato, {\it An Introduction to Infinite-Dimensional Analysis}, Springer, 2006. 
\bibitem{efk} P.I. Etingof, I.B. Frenkel, A.A. Kirillov, Jr., {\it Lectures on representation theory and Knizhnik-Zamolodchikov equations}, Providence, USA, Am. Math. Soc. (1998).
\bibitem{faddeev}L. Faddeev, {\it Modular Double of Quantum Group},  Math.Phys.Stud.21:149-156, 2000. 
\bibitem{benzvi}E. Frenkel, D. Ben-Zvi, {\it Vertex algebras and algebraic Curves}, Math. Surveys and
Monographs {\bf 88}, AMS,  2004.
\bibitem{fip}I.B. Frenkel, I.C-H. Ip, {\it Positive representations of split real quantum groups and future perspectives},  arXiv:1111.1033.  
\bibitem{zhu}I.B. Frenkel, Y. Zhu, {\it Vertex operator algebras associated to representations of affine and Virasoro algebras},  Duke Math. J. 66 (1992) 123-168.
\bibitem{ggv} I.M. Gelfand, M.I. Graev, A.M. Vershik, {\it Representations of the group $SL(2,\mathbf{R})$, where 
$\mathbf{R}$ is a ring of functions}, Russ. Math. Surv. {\bf 28} (1973) 87-132.
\bibitem{gkl} A. Gerasimov, S. Kharchev, D. Lebedev, {\it Representation theory and quantum
integrability}, Progress in Mathematics {\bf 237} (2005) 
133-156.
\bibitem{gerasimov} A. Gerasimov, S. Kharchev, D. Lebedev, S. Oblezin, {\it On a Class of Representations of Quantum Groups}, arXiv:math/0501473.
\bibitem{gmm} A. Gerasimov, A. Marshakov, A. Morozov, {\it Hamiltonian reduction of the Wess-Zumino-Witten
theory from the point of view
of bosonization}, Phys.Lett. {\bf B} 236 (1990) 269-272.
\bibitem{ivan}I.C.H. Ip, {\it The classical limit of representation theory of the quantum plane}, arXiv:1012.4145.
\bibitem{ivan2}I.C.H. Ip, {\it Positive Representations of Split Real Simply-laced Quantum Groups}, arXiv:1203.2018.
\bibitem{ivan3}I.C.H. Ip, {\it Positive representations of split real quantum groups of type $B_n$, $C_n$, $F_4$, and $G_2$ },  arXiv:1205.2940.
\bibitem{kashaev}R. M. Kashaev, A. Yu. Volkov, {\it From the Tetrahedron Equation to Universal R-Matrices},  arXiv:math/9812017
\bibitem{KL} D. Kazhdan, G. Lusztig, {\it Tensor structures arising from affine Lie algebras I,II}, J. Amer. Math. Soc. {\bf 6} (1993) 905-1011;
{\it Tensor structures arising from affine Lie algebras III,IV}, J. Amer. Math. Soc. {\bf 7} (1994) 335-453.
\bibitem{kls} S. Kharchev, D. Lebedev, M. Semenov-Tian-Shansky, {\it Unitary representation of $U_q(sl(2))$, the modular
double and multiparticle q-deformed
Toda chains}, Commun. Math. Phys. {\bf 225} (2002) 573-609.
\bibitem{quo} H.-H. Kuo, {\it Gaussian measures in Banach spaces}, Springer-Verlag, 1975. 
\bibitem{teschner} B. Ponsot, J. Teschner, {\it Clebsch-Gordan and Racah-Wigner coefficients for a continuous series of representations of $U_q(sl(2,\mathbb{R}))$}, Commun.Math.Phys. {\bf 224} (2001) 613-655.
\bibitem{teschner2} B. Ponsot, J. Teschner, {\it Liouville bootstrap via harmonic analysis on a noncompact quantum group}, hep-th/9911110.
\bibitem{ms} G. Moore, N. Seiberg, {\it Classical and quantum conformal field theory}, Commun. Math. Phys. {\bf 123} (1989) 177-254.
\bibitem{schmudgen}K. Schmudgen, {\it Operator representations of $U_q(sl(2,\mathbb{R})$}, Lett. Math. Physics {\bf 37} (1996) 211-222.
\bibitem{teschner3}J. Teschner, {\it A lecture on the Liouville vertex operators}, Int.J.Mod.Phys. A19S2 (2004) 436-458.
\bibitem{zeit} A.M. Zeitlin, {\it Unitary representations of a loop ax+b-group, Wiener measure and Gamma-function}, J. Func. Anal. {\bf 263} (2012) 529-548,  arXiv:1012.4826.
\end{thebibliography}
\end{document}